\newtheorem{theorem}{Theorem}
\newtheorem{definition}[theorem]{Definition}
\newtheorem{proposition}[theorem]{Proposition}
\newtheorem{lemma}[theorem]{Lemma}
\newtheorem{corollary}[theorem]{Corollary}
\newtheorem{remark}{Remark}

\documentclass{article}
\usepackage{amsfonts} 
\usepackage{ifpdf}
\usepackage{graphicx,amssymb}

\title{Existence of positive equilibria 
\\
for quasilinear models \\ 
of structured population}

\author{Stefano Bertoni\footnote{Dipartimento di Matematica, Universit\`a di Trento (Italy). E-mail: bertoni@science.unitn.it} 
}

\begin{document}

\maketitle 

\begin{abstract}
In this paper I prove the existence of a positive stationary solution for a generic quasilinear model of structured population. The existence is proved using Schauder's fixed point theorem. The theorem is applied to a hierarchically
size--structured po\-pu\-lation model.
\end{abstract}


\section{Introduction}
The size--structured population model IBVP (Initial Boundary Value Problem, see \cite{casa:amod}), 

in the autonomous case, has the following general form:  
\begin{equation}\label{modello1}
\cases{
	u_t + (g(x, u(t,\cdot))\,u)_x + \mu(x, u(t,\cdot))\,u=0
\cr\cr
g(0, u(t,\cdot))\,u(t,0) = \int_J \beta(x,u(t,\cdot))\,u(t,x)\,dx, 
}\end{equation}
where $x\in J=[0,\infty)$ represents \emph{age} or \emph{size}, $t\geq 0$ is \emph{time}, $u$ is the \emph{population density}, $u(t, \cdot)\in L^1(J)$ for each $t\geq 0$. 

The model equations involve the following vital rates: $\mu = \mu(x,u)$ --- mortality, $\beta = \beta(x,u)$ --- fertility and $g = g(x,u)$ --- growth rate. These coefficients depend on the size $x$ and on the total population behaviour through $u$ in a general (also nonlinear) way. 

The total population at the instant $t$ is given by  
$P(t)=\int_J u(t,x)\,dx$, the flow of the newborns is $B(t) = \int_0^\infty \beta(x, u(t, \cdot))\,u(t,x)\,dx$.
In this paper we obtain for Pbm.~(\ref{modello1}) a theorem of existence of a positive equilibrium. 

In general, however, the well-posedness of this class of PDE models is still an open question (\cite{casa:bas}, Introduction).

\par\smallskip

The first nonlinear population model  was introduced and analysed in the semi\-nal paper \cite{gmc:non} of Gurtin and MacCamy in 1974, with  nonlinearities depending only on $P(t)$. 
It was followed in the eighties by several other papers with generic nonlinearities in $u$ for the case $g=1$ e.g. by J. Pr\"uss that gave some sufficient conditions for the existence of a positive equilibrium \cite{pru:eqs,pru:ont,pru:sta}. 


%

In 2003 Diekmann et al. \cite{dgm:sap} managed  the case of nonconstant $g$ and $n$ scalar biomasses $S_1$, $S_2, \ldots S_n$ depending on $u$,  
using a very different mathematical formulation; they proved the existence of nonzero equilibria and gave bifurcation conditions. 

In 2006  Farkas e Hagen \cite{faha:sta} studied the stability of stationary solutions of the IBVP, in the case of nonlinear dependence on the total population $P$,  via linerization and semigroup and spectral methods. They give stability criteria in terms of a modified net reproduction rate. 


%
\par\smallskip

In this paper I establish Thm. \ref{teo4}, that gives sufficient conditions for the exi\-sten\-ce of a positive equilibrium for Pbm.~(\ref{modello1}), under generic dependence on $u$. 
I use a compactness hypothesis. 
I set also preliminarily some positivity and boundedness hypotheses on the coefficients $\mu$ and $g$. 
\par
\smallskip
The problem is transformed in a fixed point problem and the existence of a solution is obtained through  Schauder's fixed point theorem. 

However there is no uniqueness in general. I give a made--up counterexample. 
I give also a condition for the non--existence of positive equilibria using suitable assumptions of monotonicity on the coefficients $\mu$, $g$ and $\beta$. 
\par\smallskip
At the end of Sec.~\ref{sec:2}, I show as application the existence of a positive stationary solution for a nonlinear model of structured population of Ackleh and Ito~\cite{acit:hierc}.

In the Appendix, I resume some propositions on compactness. 

\section{Preliminaries}
\label{sec:1} 
\subsection{Notations}
$J=[0,\infty)$ is the interval of definition of $x$. 

$<\!g,f\!> =\int_J g(x)\,f(x)\,dx$ for $f\in L^1(J)$ and $g\in L^\infty(J)$. 
\par\smallskip\noindent 
$ L^1_+(J)=\{\phi\in L^1(J)\mid \phi(x)\geq 0\ {\rm a.e.}\ x\in J\}$ is the positive cone of $L^1(J)$.
\par\smallskip\noindent
Given two functions $u_1, u_2\colon [0,\infty)\to[0,\infty)$, we will write $u_1< u_2$ if  $0\leq u_1(x)\leq u_2(x)$  and $u_1\not =u_2$ a. e. $x\in J$.
The relation $<$ is a partial order on the cone 
$L^1_+(J)$.

If $e_1, e_2\in L^1(J)$ and $e_1< e_2$, then write  
$$[e_1,e_2]=\{\phi\in L^1(J)\mid e_1(x)\leq \phi(x)\leq e_2(x)\ {\rm a. e.} \ x\in J\}.
$$



Functions $f(u(\cdot))$ defined for $u\in L^1_+(J)$ will be usually briefly denoted  as $f(u)$.

\subsection{Hypotheses and definitions}
\par\medskip\noindent{\bf Hypothesis (A)}  
\begin{enumerate}

\item[a)] The functions $x\mapsto g(x,u)$, $\mu(x,u)$  
are $L^\infty(J)$ for each $u\in L^1_+(J)$ and there exist constants $\underline g, \overline g, \underline\mu, \overline\mu$: 
$$ 0<\underline g\leq g(x, u)\leq \overline g, \quad 0<\underline \mu\leq \mu(x, u)\leq \overline\mu $$ 
for each $u\in L^1_+(J)$, a. e. $x\in J$. 

\item[b)] $\beta(x, u)\geq 0$ for a.e. $x\in J$, for each $u\geq 0$ 
and there exists a constant  $\overline\beta>0$:
$\beta(x, u)\leq\overline\beta$ for each $u\geq 0$, a. e. $x\in J$.  

\item[c)] $u\mapsto g(x, u), \mu(x, u), \beta(x, u)$ are continuously depending on $u\in L^1_+(J)$ for a.e.~$x$ $\in J$.
\end{enumerate}

\par\noindent{\bf Auxiliary functions.}
For $x\in J$, $u\in L^1_+(J)$, we set: 
\begin{equation}\label{Pi}
	\Pi(x,u) := \frac{1}{g(x,u)}\,e^{-\int_0^x\frac{\mu(y , u)}{g(y , u)}\,dy }. 
\end{equation}
Under the boundedness assumptions of {\rm Hyp. (A)},  we define  the \textsl{auxiliary functions} $e_1$, $e_2$: 
\begin{eqnarray}
		e_1(x):= \frac{e^{\displaystyle -{(\overline\mu}/{\underline g})\,x}}{\overline g} 
	, \qquad 
	 e_2(x):= \frac{e^{\displaystyle -{(\underline\mu}/{\overline g})\,x}}{\underline g} . 
\end{eqnarray}

\begin{lemma}[Properties of $\Pi$]
Using the assumptions on the lower and upper bounds of $\mu$ and $g$, given in {\rm Hyp. (A)}, we obtain for each $x$, $u$: 
\begin{equation}
\label{e-diseguaglianze} 
e_1(x)  \leq \Pi(x, u)\leq   e_2(x). 
\end{equation}
Moreover $\Pi(\cdot, u)\in L^1(J)\,\cap\,L^\infty(J)$ for each $u\in L^1_+(J)$. 
\end{lemma} 
The interval 
$[e_1,e_2]$ is a closed convex subset of $L^1_+(J)$. 

\par\noindent{\bf ``Onion" set.}
\label{sec:OnionSet}
Set $\displaystyle U:=\bigcup_{\lambda>0}[\lambda\,e_1,\,\lambda\,e_2]$. 

It is simple to prove that the sets $U$ and $\overline U=U\cup \{0\}$  are convex.

%

\par\medskip\noindent{\bf Hypothesis (C)} \emph{\rm({\it Uniformly bounded variation})}. 
$$
  \forall T>0: \quad  
	\lim_{h\to 0}\sup_{u\in U}\int_0^T |g(x+h,u)-g(x,u)|\,dx = 0.  
$$
We mean that $g$ is extended as $0$ for $x<0$. 
\begin{remark}
{\rm Condition (C)} means that $\sup$ has to be considered on functions of the form $u=\lambda\,v$, with $v\in [e_1, e_2]$. Since $U\not=L^1_+(J)$ (e.g. $x^{-1/2}e^{-x}\not\in U$) this is an effective reduction of the requests. 

Under {\rm Hyp. (A)}, {\rm Condition (C)} is satisfied also for $u=0$ (therefore it holds for $u\in \overline U$) because $g(x,u)$ is continuous in $u$. 

\end{remark}

\par\medskip\noindent{\bf Hypothesis (D)} 
$$
  \forall T>0: \quad  \exists k_T>0: 
	 |g_x(x,u)|\leq k_T, \ {for\ each}\ u\in  U\ {and}\ a.e.\, x\in [0,T] .  
$$
Hyp. (D) implies Hyp. (C).

\par\medskip\noindent{\bf Hypothesis (L$_{\beta}$) ({\em Limit  of $\beta$}).} 
	For each $x\geq 0$, $$\displaystyle \lim_{\|u\|_1\to+\infty, u> 0}\beta(x, u)=0 .$$

\begin{definition}[\emph{Net reproduction function}] {\rm (Cmp. (\cite{pru:ont}, p. 330)} 
For $u\in L^1_+(J)$
\begin{equation} 
\label{gnrf}
	R(u) := \int_J \beta(x,u)\,\Pi(x, u)\,dx .
\end{equation}
\end{definition}
Under Hyp. (A), $R(u)$ is well--defined and if also (L$_{\beta}$) holds, then 
\begin{equation}
\label{R0}
	\lim_{\|u\|_1\to\infty,\, u\geq 0}R(u) = 0.
\end{equation}

\subsection{Compactness}

The (closed, convex) interval $ [e_1,e_2]\subseteq[0, e_2]\subseteq L^1_+(J)$  is invariant with respect to $\Pi$, i.e. 
$ \Pi(\cdot, [e_1,e_2])\subseteq\Pi(\cdot, [0,e_2])\subseteq [e_1,e_2]$. 
\begin{lemma}[Compactness]
\label{Pi-compact}
Under {\rm Hyp.~(A)} and {\rm(C)}, 
the function $u\mapsto\Pi(\cdot, u)$, defined on $U$ and $\overline U$ in $L^1(J)$, is compact. 
\end{lemma} 
The lemma of compactness is proved in Appendix, Sec.~\ref{sec:SectionCompactSubsetsOfE1E2}. 

\section{Existence of equilibria}
\label{sec:2}
In this section we prove the existence of a positive stationary solution $u^*$ for Pbm.~(\ref{modello1}) as fixed point of a suitable transformation of $L^1_+(J)$. 

\subsection{Stationary solutions}


The equilibria are the time--independent solutions $u=u^*(x)$ of Problem (\ref{modello1}). These are determined from   
\begin{equation}\label{pbmSs}
		\cases{
		\displaystyle 
		\frac{\partial}{\partial x}\left(g(x, u^*(\cdot))\,u^*(x) \right) + \mu(x, u^*(\cdot))\,u^*(x) =0  
		\cr \cr
		g(0, u^*(\cdot))\,u^*(0) = \int_0^\infty \beta(x, u^*(\cdot))\,u^*(x)\,dx  
		}
\end{equation}
and (see \cite{pru:ont}, Eq. (8))  they corresponds to the solutions of the functional equation 
\begin{equation}\label{u-stazionaria-1}
\label{u-stazionaria}
		u(x) = \frac{\int_0^\infty \beta(x', u)\,u(x')\,dx'}{g(x, u)}\,e^{-\int_0^x\frac{\mu(y, u)}{g(y, u)}\,dy} \quad{\rm for\ }x\in J,
\end{equation}
the only premises being $g>0$, ${\mu(\cdot, u)\over g(\cdot, u)}\in L^1_{loc}(J)$. 

This equation is translated immediately in a fixed point problem. 
\begin{proposition}\label{fixed-points}
Under {\rm Hyp.~(A)}  
the stationary solutions of  {\rm Pbm.~(\ref{modello1})} are the \emph{fixed points} of the functional ${\cal T}\colon L^1_+(J)\to L^1_+(J)$ defined as 
\begin{equation}
	\left({\cal T}\phi\right)(x) = \frac{G(\phi)}{g(x,\phi)}\,e^{-\int_0^x\frac{\mu(y, \phi(\cdot))}{g(y , \phi(\cdot))}\,dy }  
\end{equation}
and vice versa, where $G \colon L^1_+(J)\to {\mathbb R}$ is given by 
$$G(\phi) = \int_J \beta(x, \phi(\cdot))\,\phi(x)\,dx, $$
for $\phi\in L^1_+(J)$. 
\end{proposition}
The functional equation  $u = {\cal T}\,u$ 
can be written in a more compact form as 
\begin{equation}
\label{forma1}
	u(x) = G(u(\cdot))\,\Pi(x, u(\cdot)),  
\end{equation}
that we discuss. 

\begin{theorem}[Existence of equilibria]
\label{teo4}
Assume {\rm Hyp.~(A)} and {\rm (C)}. 
Suppose there is a constant $\rho_0>0$ such that for $u\in L^1_+(J)$, $\|u\|_1\geq \rho_0$ implies $R(u) \leq 1$. 
\\ 
If $R(0)>1 $ 
then {\rm Problem (\ref{modello1})} admits at least a positive stationary solution.  

The solution satisfies the functional equation 
\begin{equation}
	u^*(x) = \lambda^*\,\Pi(x,u^*(\cdot)), 
\end{equation}
where $\lambda^*>0$  is a suitable number and 
the corresponding population is constant and given by $\displaystyle P^*=\lambda^*\,\|\Pi(\cdot, u^*)\|_1$.  
\end{theorem}
From (\ref{R0}) we have the following statement:
\begin{corollary}
\label{teo4b}
Under {\rm Hyp.~(A), (C)} and {\rm(L$_\beta$)},  
if $R(0)>1 $ then {\rm Problem (\ref{modello1})} admits a positive stationary solution.  

\end{corollary}
\subsection{Proof of Thm.~\ref{teo4}}
Prop.~\ref{fixed-points} reduces the search for equilibria of Pbm.~(\ref{modello1}) to Eq.~(\ref{forma1}).  
\par\smallskip 
$G(0)=0$ gives the trivial equilibrium $u = 0$ so we exclude this case. 


The proof is divided into two steps. 
\par\smallskip 
{\bf{\rm(i)} Splitting variables.} 
Consider Eq. (\ref{forma1}): 
 assume that $u$ is a solution of $u=G(u)\,\Pi(\cdot, u)$. 
\\
Set $\lambda := G(u)\ (\not= 0)$ and $\displaystyle v = \frac{1}{\lambda}\,u$. 
\\ 
By substitution we obtain: $\lambda\,v =\lambda \Pi(x, \lambda\,v)$ and $\lambda = \int_0^\infty \beta(x, \lambda\,v)\, \lambda\,\Pi(x, \lambda\,v)\,dx$ so that $1 = \int_0^\infty \beta(x, \lambda\,v)\, \Pi(x, \lambda\,v)\,dx$. Hence $(v, \lambda)\in  [e_1,e_2]\times(0,\infty)$ is a solution of the system: 
\begin{equation}
\label{sis2}\label{sis3}
	\cases{
	v(x)=\Pi(x, \lambda\,v(\cdot)),  
	\cr\cr 	  
	R(\lambda\,v)=1.}
\end{equation}
Vice versa, if $(v, \lambda )$ is a solution of (\ref{sis2}), then $u=\lambda\,v$ is a solution of the equation $u=G(u)\,\Pi(\cdot,u)$.

The condition $R(0)>1$   implies that $\lambda^*\not=0$. For each solution $(v,\lambda)$ of Pbm.~(\ref{sis2}) we have $(v, \lambda)\in  [e_1,e_2]\times(0,\infty)$. 


{\bf{\rm(ii)} Fixed point.}
In this step we apply Schauder's fixed point theorem --- see \cite{evans:pde}, \cite{zei:nfa}. We write Pbm.~(\ref{sis2}) in the form  
\begin{equation}
\label{system-4}
\cases{v(\cdot) = \Pi(\cdot, \lambda\,v),\quad v\in [e_1, e_2],
\cr 
\lambda = \max\{\lambda + R(\lambda\,v) - 1 ;\, 0\},\quad \lambda\geq 0} 
\end{equation}
that is $ (v, \lambda)= A\Big((v,\lambda)\Big)$, with $(v,\lambda)\in [e_1,e_2]\times(0,\infty)$ and $A$ defined by the second members of (\ref{system-4}).


The map $u\mapsto\Pi(\cdot, u)$ is continuous and compact on $U$; the function $R(u)$ is continuous and bounded from $L^1_+(J)$ to $(0, \infty)$, since $ 0< R(u)\leq \overline \beta\,\|e_2\|_1$; therefore $A\colon [e_1, e_2]\times (0, \infty)\to L^1(J)\times (0, \infty)$ is continuous and compact. 

$A_1 (v,\lambda):= \Pi(\cdot, \lambda\,v)$ has image in $[e_1,e_2]$.  
\\ 
Now prove that for a fixed $\displaystyle M > \frac{\rho_0}{\|e_1\|_1}$, $\displaystyle \frac{\rho_0}{\|e_1\|_1} + \overline\beta\,\|e_2\|_1 - 1$, then \\ 
 $A_2(v,\, \lambda):= \max\{ \lambda + R(\lambda\,v)-1;\, 0 \}$ 
maps $[e_1,e_2]\times [0,M]$ on $[0,M]$. 

If $\displaystyle\frac{\rho_0}{\|e_1\|_1}\leq \lambda \leq M$, then $\displaystyle\lambda \geq \frac{\rho_0}{\|v\|_1} $ and $R(\lambda\,v )\leq 1$, 
so that
$$\lambda + R(\lambda\,v) - 1 \leq 1 + \lambda-1 = \lambda \leq M. $$

If $0\leq \lambda < \rho_0/\|e_1\|_1$, then 
$\displaystyle\lambda+R(\lambda\,v)-1 \leq \frac{\rho_0}{\|e_1\|_1} + \overline \beta\,\|e_2\|_1 - 1 \leq M$.  

So $A$ maps $[e_1, e_2]\times [0,M]$, a closed convex subset of $L^1(J)\times (0,\infty)$, in itself. 
\par\medskip\noindent 
Since $A$ is compact,  
by Schauder's fixed point theorem, Eq.~(\ref{system-4}) has at least a fixed point $(v^*, \lambda^*)\in [e_1 e_2]\times [0,M]$ and it is different from $0$ for the initial remark; $(v^*, \lambda^*)$ is a fixed point also for Eq.~(\ref{sis2}).  
\par\medskip 
\noindent 
Finally, Eq. (\ref{forma1}) is satisfied by $u^* = \lambda^*\,v^*$ and the corresponding stationary population is $$P^*=\int_J u(x)\,dx = \lambda^*\int_J v^*(x)\,dx .$$

\begin{remark}
$R(0)>1$ implies $\overline \beta\,\|e_2\|_1> 1$, therefore in the proof it is possible to assume 
$M=\frac{\rho_0}{\|e_1\|_1} + \overline \beta\,\|e_2\|_1 - 1$ and to have the estimate $P^* \leq M\,\|e_2\|_1$. 
\end{remark}

\subsection{A counterexample}
{\rm Thm.~\ref{teo4}} is a sufficient condition but not a ne\-ces\-sary one. We can have also $R(0)<1$ if there exists $u_0\in L^1_+(J)$ such that $R(u_0)>1$. In this case it is possible to need other conditions on $u_0$ to prove a statement of existence. 
The idea is to construct explicitly an example with a positive equilibrium but $R(0) <1$. 

Set  $\mu(x,u)=g(x,u)=g$ so that $\displaystyle \Pi(x,u) = \frac{1}{g}\,e^{-x} $, independent of $u$. 

Define $e_0(x):= e^{-x}$. 
Take $F\colon L^1_+(J)\to {\mathbb R}, u\mapsto F(u)$, such that $F(0)<1$,  $F(e_0)=1$ and $\displaystyle\lim_{\|u\|_1\to\infty,\ u>0} F(u)=0$, $F$ continuous but obviously nonmonotonic. 

Now set $\beta (x,u)=2\,g\,(1-e^{-x})\,F(u)$ so that $R(u) = F(u)$.

Then $R(e_0)=1$ and $u = e_0$ is a solution of the fixed point equation and a positive equilibrium. 
\par\smallskip\noindent 
As example of function $F$ we can take $F(u):= f(\|u\|_1)$. where 
\begin{equation}
	f(a) := \cases{\displaystyle\frac{1}{2}+ 3a & for $0\leq a\leq\frac{1}{2}$, 
	\cr
	\displaystyle
	3-2\,a & for $ \frac{1}{2}<a\leq \frac{5}{4}$,  
	\cr
	\displaystyle
	\frac{e^{5/4}}{2}\,e^{-a}& for $a> \frac{5}{4}$. 
	\cr} 
\end{equation}
In this case we have \emph{two} positive equilibra, 
$u(x)=e^{-x}$ and $u(x)=\frac{1}{6}\,e^{-x}$, corresponding to the two solutions of $f(a)=1$, i.~e. $a=1$, $a=1/6$.  

\subsection{A nonexistence result and a sufficient and necessary condition}
\label{sec:MonotoneCase}

Under suitable monotonicity hypotheses, $R(0)> 1$ becomes a necessary and sufficient condition.

%
A function $f$, defined on ordered spaces, is \emph{increasing} if $u_1<u_2$ implies $f(u_1)<f(u_2)$. The other monotonicity definitions are extended in the same ways.

Now assume $u\in L^1_+(J)$ in the following statements. 
\par\medskip\noindent{\bf Assumption (M) (\emph{Monotonicity}) }
\begin{itemize}
	\item $u\mapsto \mu(x, u)/g(x, u) $ is nondecreasing (or increasing) for each $x\geq 0$ 	(\emph{morta\-li\-ty--growth ratio}),
	\item $u\mapsto \beta(x, u)/\mu(x, u)$ is decreasing (or nonincreasing) for each $x\geq 0$ (\emph{fertili\-ty--mor\-ta\-li\-ty ratio}), 
	\item $x\mapsto \beta(x, u)/\mu(x,u)$ is nondecreasing (or increasing) for each $u$.
\end{itemize} 
The hypotheses between parentheses are in alternative: $u\mapsto \beta/\mu$ must be strictly decreasing and the other two functions are only nondecreasing, or, vice versa, $u\mapsto \beta/\mu$ nonincreasing and the others have to be two strictly increasing. 

To prove the nonexistence condition we need the following statement:  
\begin{lemma}[Monotonicity] 
\label{mon-u}
Assume {\rm Hypotheses~(A)}, {\rm (C)}, {\rm(L$_{\beta}$)} and \emph{As\-sum\-ption (M)}. 

Then the functional $R\colon L^1_+(J)\to (0,\infty)$ is continuous, decreasing and 
$$\lim_{\|u\|_1\to+\infty, u> 0} R(u) = 0.$$
\end{lemma}
I do not give the details of the proof of this lemma, but the main idea is to write $R(u)$ as $\int_J dx\, \frac{\beta(x,u)}{\mu(x,u)}\frac{\mu(x,u)}{g(x,u)}\,e^{-\int_0^x \frac{\mu(y,u)}{g(y,u)}\,dy}$ and to study the properties of monotonicity of the integral $\int_J dx\, h(x)\,f(x)\,e^{-\int_0^x f(y)\,dy}$ with respect to suitable $f$ and $h$.

For a detailed proof, see Bertoni \cite{ber:monotone}.

\begin{proposition}[Non existence of positive stationary solutions] 
\label{nonexist} $\,$\\ 
Under Hypotheses of {\rm Lemma~\ref{mon-u}}, 
if $R(0)\leq 1$  
then {\rm Pbm.~(\ref{modello1})} has no positive stationary solutions.
\end{proposition}
{\it Proof.}
If $R(0)\leq 1$ then $R(u)=1$ does not have positive solutions by monotonicity. 

Since existence of positive equilibria is equivalent to positive solutions of $u=G(u)\,\Pi(\cdot, u)$ and so of Eq.~(\ref{sis3}), the conclusion follows. 
\par
\medskip\noindent 
As consequence, Condition $R(0)> 1 $ becomes a \emph{necessary and sufficient condition of existence of positive equilibria}  for Pbm.~(\ref{modello1}) under Hyp.~(A), {\rm (C)}, {\rm(L$_{\beta}$)} and (M). 

\subsection{Applications}
\label{sec:Applications} 
Ackleh e Ito~\cite{acit:hierc} consider a hierarchically
size-structured population model that can be reported to Eq.~(\ref{modello1}). They proved existence of measure-valued solutions for the  Cauchy problem. We give a condition of existence of a stationary positive solution for a simple case of this model, by taking  
\begin{equation}
\label{hierc-g}
	g(x, u(\cdot))= \underline g + (\overline g - \underline g)\,e^{-\int_x^\infty u(y)\,dy}. 
\end{equation}
Hyp.~(D) is equivalent to
 $$ \forall T>0: \quad  
\sup_{u\in U} {\rm\,\ ess\!\!}\sup_{0\leq x \leq T}  |g_x(x, u)|<\infty  
$$
that is, for (\ref{hierc-g}):   
\begin{equation}
	\forall T>0: \quad  
	\sup_{u\in U} {\rm\,\ ess\!\!}\sup_{0\leq x\leq T}| e^{-\int_{x}^\infty u(y)\,dy}\cdot u(x)|_\infty<\infty. 
\end{equation}
For $u=\lambda\,v$ with $v\in[e_1,e_2]$ we use the inequality  
$\displaystyle\sup_{\lambda>0}\lambda\,e^{-\alpha\,\lambda}=\frac{1}{\alpha\,e}$: therefore  
$$ 
e^{-\int_{x}^\infty u(y)\,dy}\cdot u(x) = \lambda\,v(x)\,e^{-\lambda\int_{x}^\infty v(y)\,dy}\leq \frac{v(x)}{e\,\int_{x}^\infty v(y)\,dy}\leq \frac{e_2(x)}{e\int_T^\infty e_1(y)\,dy}<\infty.$$

Assume $\mu$ and $\beta$ to satisfy Hyp. (A) and (L$_\beta$). The other conditions on $g$ of Cor.~\ref{teo4b} are trivially satisfied, so we obtain the existence of at least one positive stationary solution if 
$$
\int_J dx\, \beta(x,0)\, e^{-\int_0^x \frac{\mu(y, 0)}{g(y,0)}\,dy}>\overline g. $$ 

\appendix
\section*{Appendix}

\section{Compactness conditions} 
 
 
As well known, the conditions for the relative compactness of a set $W$ in $L^1(0,\infty)$ are given by the Riesz--Kolmogorov Theorem. We use the following version:  
\par\medskip\noindent 
\emph{i)} $W $ is bounded; 
\\
\emph{ii)} 
$\displaystyle\lim_{T\to \infty}\sup_{w\in W}\int_{x>T} |w(x)|\,dx  = 0$. 
\\
\emph{iii)} 
$
\displaystyle\lim_{h\to 0}\sup_{u\in W}\int_{0}^T |w(x+h)-w(x)|\,dx =0
$ for each $T>0$.

\noindent 
Sets of continuous, uniformly bounded variation functions in $L^1(0,\infty)$ are (relatively) compact. 

\section{Compactness of $\Pi$ {\rm (Proof of Lemma~\ref{Pi-compact})}}
\label{sec:SectionCompactSubsetsOfE1E2}

For each $u\in L^1_+(J)$, the function $\Pi(x,u)$ defined by (\ref{Pi}) has the following properties: 
\begin{enumerate}
	\item $\Pi(\cdot,u)\in [e_1, e_2]$, that implies  \emph{(i)} and \emph{(ii)} of the Riesz--Kolmogorov Theorem; 
	\item $x\mapsto \Pi(x,u) $ is continuous. 
\end{enumerate}
Now we prove \emph{(iii)} for $u\in \overline U$. Let be $T>0$, $h>0$:  
\begin{eqnarray*}
	&& \int_0^T |\Pi(x+h,u) - \Pi(x,u)|\,dx \leq 
	\\ 
	&& \leq \int_0^T dx\, \frac{e^{-\int_0^{x+h} \frac{\mu(y, u)}{g(y,u)}\,dy}}{g(x+h, u)}\left({e^{\int_x^{x+h} \frac{\mu(y, u)}{g(y,u)}\,dy}} - 1\right) + \\ &&
	+ \int_0^T dx\,\frac{e^{-\int_0^x \frac{\mu(y, u)}{g(y,u)}\,dy} }{g(x,u)\,g(x+h,u)}\,|g(x+h,u)-g(x,u)|
	\leq
	\\ && \leq 
	\frac{T\,\overline\mu}{\underline g^2} \,h + \frac{T}{\underline g^2}\,\int_0^T dx\,|g(x+h,u)-g(x,u)|, 
\end{eqnarray*}
therefore, using Hyp.~(C) for $u\in \overline U$, this completes the proof. The case $h<0$ is managed analogously. 

We obtain the $\Pi$  sends $ U $ in a relatively compact subset of $[e_1,e_2]$ in the norm of $L^1(J)$ i. e. the set $\Pi(\cdot,U)$ is relatively compact.

\thanks 
\paragraph{\it Acknowledgements.}
\label{sec:Acknowledgements}
I thank A. Pugliese (Univ. of Trento) for the discussions and the remarks, the example in Sec. \ref{sec:Applications} and a suggestion to simplify the proof of Thm.~\ref{teo4}.


\begin{thebibliography}{24} 

\bibitem{ber:monotone} S. Bertoni, \emph{Monotonicity of a Class of Integral Functionals}, 
{\it preprint} arXiv:1503.05390 (2015)

\bibitem{acit:hierc} A. S. Ackleh, Kazufumi Ito, \emph{Measure--valued solutions for a hierarchically size--structured population}, J. Differential Equations {\bf 217} (2005), 431--455.

\bibitem{casa:amod} A. Calsina, Joan Salda\~{n}a, \emph{A model of physiologically structured population dynamics with a nonlinear individual growth rate}, J. Math. Biol. {\bf 33} (1995), 335--364.
\bibitem{casa:bas} --, \emph{Basic Theory for a Class of Models of Hierarchically Structured Population Dynamics with Distributed States in the Recruitment}, Mathematical Models and Methods in Applied Sciences {\bf 16}, No. 10 (2006), 1695--1722.
%
%
%
\bibitem{dgm:sap} O. Diekmann, M. Gyllenberg and J.A.J. Metz, {\it Steady-state analysis of structured population models}, Theoretical Population Biology {\bf 63} (2003), 309-–338. 
\bibitem{evans:pde} L. C. Evans, ``Partial Differential Equations'', American Mathematical Society, Providence, 1998. 
\bibitem{faha:sta} Jozsef Z. Farkas, Thomas Hagen, \emph{Stability and regularity results for a size-structured population model}   J. Math. Anal. Appl. {\bf 328} (2007), 119--136. 
\bibitem{gmc:non} M. E. Gurtin, R. C. MacCamy, \emph{Nonlinear age--dependent population dynamics}, Arch. Rat. Mech. Anal. \textbf{54} (1974), 281--300. 


\bibitem{pru:eqs} Jan Pr\"uss, \emph{\it Equilibrium solutions of age-specific population dynamics of several species},  J. Math. Biol. {\bf 11}, 1 (1981), 	65--84.

\bibitem{pru:ont} --, \emph{On the qualitative behaviour of populations with age-specific interactions}, Comp. Math. Appls. {\bf 9}, 3 (1983), 327--339 

\bibitem{pru:sta} --, \emph{Stability analysis for equilibria in age-specific population dynamics}, Nonlin. Analysis TMA {\bf 7}, 12 (1983), 1291--1313. 


\bibitem{zei:nfa} E. Zeidler, ``Nonlinear Functional Analysis and its Applications, Vol. I -- Fixed Point Theorems'', Springer--Verlag, New York, 1986. 

\end{thebibliography}
\end{document}